\documentclass{amsart}
\usepackage{hyperref}

\usepackage{verbatim}

\newcommand{\bio}[2]{\noindent{\bf #1.} \emph{#2}\smallskip}
\newcommand{\quot}[2]{\noindent{\bf #1:} #2\smallskip}

\title{An oral history of FPSAC}
\date{October 13, 2026}
\author{Sam Hopkins}
\address{Department of Mathematics, Howard University, Washington, DC 20059, USA}
\email{samuelfhopkins@gmail.com}

\begin{document}

\begin{abstract}
As its 40th anniversary nears, we present a history of the International Conference on Formal Power Series and Algebraic Combinatorics (FPSAC), through interviews with figures from its past, present, and future.
\end{abstract}

\maketitle

The International Conference on Formal Power Series and Algebraic Combinatorics (FPSAC) is an annual mathematical conference in the area of enumerative and algebraic combinatorics. The first FPSAC conference was held in 1988 in Lille, France, the second in 1990 in Paris, and there has been a FPSAC conference in some city around the world every year since then.\footnote{The 2020 conference was held online because of the COVID pandemic.} FPSAC has taken place on five continents in 19 different countries: France, Canada, Italy, the United States, Austria, Spain, Russia, Australia, Sweden, China, Chile, Iceland, Japan, South Korea, the United Kingdom, Slovenia, Israel, India, and Germany. The 38th FPSAC conference took place in Seattle, Washington, USA in July of 2026, and the 39th will take place in Galway, Ireland in July of 2027.

In many ways, FPSAC grew up together with the field of enumerative and algebraic combinatorics. Enumerative and algebraic combinatorics really emerged as an independent discipline in the 1960s and 70s -- see for instance Richard Stanley's lively account in~\cite{stanley2021enumerative}. Since the 1980s, enumerative and algebraic combinatorialists have uncovered deep connections between their research and other areas of pure math, like commutative algebra, algebraic geometry, representation theory, probability, and complex analysis, as well as other scientific fields, like computer science, statistical mechanics, other parts of theoretical physics, and mathematical biology. Recent prestigious awards, such as the 2022 Fields Medal awarded to June Huh, and the 2022 Steele Prize for Lifetime Achievement awarded to Richard Stanley, testify to the increasing recognition of enumerative and algebraic combinatorics in the mathematical community at large. (June Huh was an invited speaker at FPSAC in 2014; Richard Stanley has been an invited speaker on three occasions.)

Nowadays, hundreds of mathematicians gather every summer for FPSAC, many of them of a younger generation who may not know much about the history of the conference. The purpose of this article is to record that history. For example, we will explain the ``formal power series'' part of the FPSAC name. We will also discuss the conference's francophone roots. And we will touch on the present and future of the conference series as well. Since the best way to learn history is from those who experienced it, our method will be to interview several figures from the history of FPSAC: Christophe Reutenauer, Maylis Delest, Gilbert Labelle, Fran\c{c}ois Bergeron, H\'{e}l\`{e}ne Barcelo, Nantel Bergeron, Christian Krattenthaler, Mike Zabrocki, Bridget Tenner, and Jang Soo Kim. These mathematicians shaped the past, present, and future of FPSAC. So let's hear from them. 

\bio{Christophe Reutenauer}{Christophe Reutenauer is a professor in the Department of Mathematics at the Universit\'{e} du Qu\'{e}bec \`{a} Montr\'{e}al (UQAM). A former student of the \'{E}cole Normale Sup\'{e}rieure (Paris), he received his PhD in math from Universit\'{e} Pierre-et-Marie-Curie (UMPC) -- Paris VI in 1980 under the supervision of Jean Berstel. He has made fundamental contributions to many parts of algebraic combinatorics, including symmetric and quasi-symmetric functions, Hopf algebras, combinatorics of words, and free Lie algebras. He was a co-organizer of the first FPSAC conference in Lille, France in 1988.}

\quot{S.H.}{Along with G\'{e}rard Jacob, you put together the first FPSAC/SFCA conference in 1988 at the University of Lille. But this first conference did not have the current name of ``Formal Power Series and Algebraic Combinatorics'' (``S\'{e}ries Formelles et Combinatoire Alg\'{e}brique''). Rather, it was a conference about ``Algebraic and Computing Treatment of Noncommutative Power Series'' (``Traitements Alg\'{e}briques et Informatiques des S\'{e}ries Formelles Non Commutatives''). What motivated you to organize a conference on this specific topic at that time?}

\quot{C.R.}{G\'{e}rard Jacob and me, following Marcel-Paul Sch\"{u}tzenberger and Michel Fliess, worked on noncommutative formal power series. Fliess applied these series to control theory, where they are often called Chen--Fliess series.  All this was the motivation of the first conference, principally organized by Jacob (I am a lazy organizer), and justifies its title. Jacob was professor at Lille at this time (but living in Paris, he was a so-called ``turbo-prof'') and I was researcher at CNRS in Paris (my lab was the LITP), about to leave for Montr\'{e}al, UQAM. The next conference in Paris in 1990 had also this title and flavor. G\'{e}rard Jacob died in 2024.}

\quot{S.H}{When you organized that first conference in Lille, did you envision it continuing as an annual series for the next 40 years?}

\quot{C.R.}{Certainly not.}

\quot{S.H.}{How did FPSAC evolve from that first conference about noncommutative power series to a much broader conference series focused on all aspects of algebraic combinatorics?}

\quot{C.R.}{I wouldn't say it became ``much broader,'' but rather experienced a change of direction, towards algebraic combinatorics. The noncommutative aspects of formal power series were somewhat forgotten. I think the title changed for the Bordeaux conference in 1991.}

\quot{S.H.}{The participants of the early FPSAC conferences were mostly French speakers, from both France and Canada. Do you remember for that first FPSAC conference in Lille, were the talks given in English or in French, or a mixture?}

\quot{C.R.}{A mixture; I would say the abstracts are written in the language of the talk.}

\quot{S.H}{And when did FPSAC start to move away from its francophone roots?}

\quot{C.R.}{The roots are francophone, and even French: all speakers, but one who is Italian, in Lille 1988 are French, as you can see on the abstracts. But the organization became fully international for the Montr\'{e}al conference in 1992.}

\quot{S.H.}{Do you have any favorite memories from those early FPSAC years?}

\quot{C.R.}{I was invited speaker in 1994 (Rutgers), co-chair of the program in 1992 (Montr\'{e}al), and I remember the wine castle in 1991 (Bordeaux). In 1992, we had the luck to have as invited speaker I.M.~Gelfand, one of the most famous Russian mathematicians of the 20th century; also, Louis Solomon was there, and I was happy to meet him. I remember also Florence 1993, because it is Florence.}

\bio{Maylis Delest}{Maylis Delest is a retired professor in the Laboratoire Bordelais de Recherche en Informatique (LaBRI) at the Universit\'{e} de Bordeaux. She received her th\`{e}se d'\'{e}tat in math from the Universit\'{e} de Bordeaux in 1987 under the supervision of Xavier G.~Viennot. Her research is in enumerative combinatorics and computer science, with a focus since the year 2000 on the visualization of information. She presented at the first two FPSACs, and then organized the 3rd FPSAC in Bordeaux, France in 1991. She was a founding member of the FPSAC executive committee.}

\quot{S.H.}{I know you attended the first two FPSAC conferences, and organized the 3rd one in Bordeaux in 1991. How did you get involved with FPSAC at its inception?}

\quot{M.D.}{First, it is important to understand the context: in France at the time, almost all combinatorics researchers worked in computer science laboratories. In Lille, the conference was organized by G\'{e}rard Jacob, who worked at LIFL, and Christophe Reutenauer, now a professor at UQAM, who worked at LIAFA. The ``I'' in both acronyms stands for ``Informatique.'' Most researchers also had a master's degree in computer science, with training split roughly equally between mathematics and computer science, and many completed PhDs in computer science as well. In 1988 and 1990, about 30 people attended the meetings. In September 1990, Daniel Krob visited Bordeaux, and we discussed organizing the next conference.}

\quot{S.H.}{This third FPSAC in Bordeaux was the first to use the name ``Formal Power Series and Algebraic Combinatorics,'' and moved beyond an exclusive focus on power series to include other aspects of enumerative and algebraic combinatorics as well. How was the decision made to broaden the scientific scope of the conference?}

\quot{M.D.}{Jean Berstel was in Bordeaux at the end of 1990 and felt it was a good idea to broaden the conference. I should point out that I was not the chair of the program committee -- G.~Jacob and P.~Leroux were. To define the wider scope, we brought together our complementary expertise: mine in formal computation and enumerative combinatorics, Daniel Krob's in algebraic combinatorics, G\'{e}rard Jacob's in symbolic computation, and Jean Berstel's broad knowledge of combinatorics. This collaboration shaped both the call for papers and the conference title. We also chose to organize FPSAC like a computer science conference, with full-paper submissions reviewed by a program committee, although that format has since been abandoned.}

\quot{S.H.}{What was organizing that conference like? Do you have any fun memories?}

\quot{M.D.}{Here were some fun elements. \\ \indent We expected 50 persons, and 130 attended the conference, surely due to some senior folks like G.~Viennot and R.~Cori and also the names of the members of the program committee. (In Firenze in 1993, organized by Renzo Pinzani, the number of participants was more than 200.) \\ \indent We got so many submission that one of my colleagues, Jean-Guy Penaud, organized the first GASCom conference. (This conference is still existing.) \\ \indent The conference began on May 2, but May 1st was a public holiday: there were no taxis or buses, and many restaurants were closed. With my colleagues, we drove participants to their hotels and hosted some of them at our homes; I personally welcomed about 20 persons with a fantastic duck dinner (unfortunately for Julian West who was vegan). \\ \indent We had long discussions on ``internet and publication'' with Herbert Wilf and he decided to launch an \emph{Electronic Journal of Combinatorics}. We were discussing frequently about the technical points by mail during the set up.}

\quot{S.H.}{How did the FPSAC executive committee start? How long did you serve?}

\quot{M.D.}{After several years we decided to create an executive committee with Daniel Krob and Pierre Leroux. I think Jean Berstel was also there; I remember that we were five. I left it surely in 2006 when I was not involved any more in combinatorics.}

\quot{S.H.}{Your own research interests have evolved over the course of your career, from enumerative combinatorics to computer science. Can you talk a little about that?}

\quot{M.D.}{My field was enumerative combinatorics. At one point I needed to draw the objects I was counting and combine that with symbolic computation. That's how, with PhD students, I developed the CalICo software. I don't remember why, but I met Ivan Herman (CWI) and we talked about an emerging discipline: information visualization. Gradually, I found myself back in the world of computer science. At the same time,  Jean-Daniel Fekete was developing a research project on information visualization at Orsay, near Paris. Thus, I continued my career in this field. If you want to explore the link with combinatorics, look at my slides ``From polyominoes to information visualization''\footnote{\url{https://www.labri.fr/perso/bousquet/XGV2012/MaylisDelest.pdf}} (for a journey in honor of Xavier Viennot) or my paper ``Efficient drawing of RNA secondary structure'' with Auber, Domenger, and Dulucq (linking infovis, combinatorics, biology).}

\quot{S.H.}{Do you think FPSAC does a good job at connecting enumerative and algebraic combinatorics with other related fields like computer science? When you were on the executive committee, was that an intentional focus?}

\quot{M.D.}{Yes and yes. We took care to maintain a balanced representation of the different research themes in the conference program. This is why we introduced a poster session, allowing authors of papers that could not be accepted -- sometimes simply because of the need to maintain this balance -- to present their work.}

\quot{S.H.}{When was the last FPSAC conference you attended? Have you kept up with more recent developments in enumerative and algebraic combinatorics?}

\quot{M.D.}{Taormina, 2005, in honor of my friend Adriano Garsia!\footnote{Adriano Garsia died in 2024; see~\cite{garsia2026memories} for remembrances of him.} I meet him the first time in 1983. When I first visited him in San Diego, we began writing a cookbook together, ``Recettes des quatre continents,'' and had a lot of mathematics questions on the beach. As for recent developments, I have often been asked to referee papers on polyominoes. I am retired now.} 

\bio{Gilbert Labelle}{Gilbert Labellle is a professor emeritus in the Department of Mathematics of the Universit\'{e} du Qu\'{e}bec \`{a} Montr\'{e}al (UQAM). He received his PhD in math from the Universit\'{e} de Montr\'{e}al in 1969 under the supervision of Qazi Ibadur Rahman. His research interests evolved during his career from analysis to formal power series and enumerative combinatorics, including the theory of species. He was an invited speaker at the 3rd FPSAC in Bordeaux in 1991.}

\quot{S.H.}{Before we talk about FPSAC, let me ask about your research. You started out in analysis, then moved into combinatorics -- how did that transition happen?}

\quot{G.L.}{From polynomials to convergent power series in analysis, I became interested in formal power series. Moreover since the theory of species provides structural aspects to various kinds of formal power series, I became immediately drawn to enumerative and algebraic combinatorics as one of my main research interests.}

\quot{S.H.}{You were an invited speaker at the FPSAC conference in Bordeaux in 1991. How did you come to be invited to speak there, and what did you talk about?}

\quot{G.L.}{The theory of combinatorial species was introduced by Andr\'{e}
Joyal in 1981. It takes into account the whole automorphism groups of combinatorial structures, not only their enumeration. In particular, this theory interconnects P\'{o}lya theory and various kinds of formal power series. I became immediately interested in this theory and I used part of it to give a combinatorial proof of the Lagrange inversion formula. Later, with Pierre Leroux, we organized an international conference on enumerative combinatorics, in which combinatorial species had an important role, in 1985 in Montr\'{e}al. Following that conference, Pierre Leroux became one of the promoters of FPSAC conferences. (Pierre Leroux died in 2008; I wrote about his mathematical life in~\cite{labelle2021leroux}.) The organizing committee of the 1991 FPSAC in Bordeaux then invited me to give a lecture with the following title, in which automorphism groups of structures and species theory are central: ``Sur la sym\'{e}trie et l'asym\'{e}trie des structures combinatoires.''}

\quot{S.H.}{Which other FPSAC conferences did you participate in after that?}

\quot{G.L.}{I participated in 15 FPSACs from Bordeaux 1991 to Paris 2013, with over 20 abstracts submitted during this period.}

\quot{S.H.}{A great conference experience is when chatting with another participant leads to a new collaboration. Do you have any collaborations that started at FPSAC?}

\quot{G.L.}{Yes, I believe two such research collaborations originating from a FPSAC conference were my paper with Ira Gessel, ``Lagrange inversion for species'' and my paper with Leroux, Pergola, and Pinzani, ``Stirling numbers interpolation using permutations with forbidden subsequences.''}

\quot{S.H.}{Did you have a favorite FPSAC conference?}

\quot{G.L.}{I much enjoyed FPSAC '07, Nankai University, Tianjin, China. The organization was fabulous!}

\bio{Fran\c{c}ois Bergeron}{Fran\c{c}ois Bergeron is the Chair of the Department of Mathematics of the Universit\'{e} du Qu\'{e}bec \`{a} Montr\'{e}al (UQAM). He received his PhD in math from the Universit\'{e} de Montr\'{e}al in 1986 under the joint supervision of Pierre Berthiaume and Pierre Leroux. He has made important contributions to several areas of algebraic combinatorics, including the theory of species, descent algebras, Macdonald polynomials, and rational Catalan combinatorics. Together with his brother and collaborator Nantel Bergeron, he participated in the 3rd FPSAC conference in Bordeaux, France in 1991, and has attended many FPSACs since.}

\quot{S.H.}{The first FPSAC conference you attended was the one in Bordeaux in 1991. How did you learn about this conference, and what motivated you to participate?}

\quot{F.B.}{There was an ongoing series of exchanges between our gang in Montr\'{e}al, the combinatorial groups in France--Germany (Paris, Strasbourg--Erlangen and Bordeaux), as well as the MIT gang. Our group in Montr\'{e}al was involved in the program committee. It was thus natural for me to participate. In a sense, that meeting is the true birth of FPSAC as an intentional series of conferences. A precursor was the meeting that was organized in Montr\'{e}al in 1985. If you look up the proceedings (see~\cite{montreal1985proceedings}), you will see many links with FPSAC. The proceedings include a list of participants that you will recognize among those of FPSAC. Among the speakers, besides me, there were: Bressoud, Viennot, Foata, Garsia, Gessel, Joyal, Kerber, Reutenauer, Stanley, and Strehl. Maybe this was one of the first meetings really dedicated to algebraic combinatorics. This all made natural to participate in the Bordeaux FPSAC in 1991.}

\quot{S.H.}{How many FPSAC conferences would you estimate you've attended?}

\quot{F.B.}{If I did not forget one, I count 16 at least instances of actual participation, with me giving regular talks on six occasions, and twice more being invited.}

\quot{S.H.}{Were you ever involved with the organization of a FPSAC conference? I know that your brother Nantel organized the conference in Toronto in 1998, and was on the executive committee for some years.}

\quot{F.B.}{I do not remember being an organizer. I was rather 6 times on the program committee (3 of those without attending), and I was chair of the program committee on another occasion (which I attended). I would typically not be involved with one of the committees if I gave a talk.}

\quot{S.H.}{You have been at the center of many important developments in enumerative and algebraic combinatorics over the past 40 years. How would you say that FPSAC has contributed to these kinds of developments in this field?}

\quot{F.B.}{Among all conferences that are involved in the development of enumerative and algebraic combinatorics, FPSAC is by far the one with the most impact. Just as an indication, the list of invited speakers contains a significant portion of the current international leaders in the field, which is entirely augmented by the list of contributed talks as well as posters. One clearly sees the emergence of new leaders along the years.}

\quot{S.H.}{Do you have a favorite memory from those early FPSAC years?}

\quot{F.B.}{There are lots of fun stories. One that I remember well is being out of breath and overfull after the lunch break in Firenze, just before my invited talk. The local organizers had invited me to a very nice restaurant down the hill (probably with some wine in the mix). Being late, we had to run up for my talk. As you can imagine, this was not an ideal way for me to give a talk. In fact, this may have been one of my worst talks, but the memory is still nice.}

\quot{S.H.}{In your view, what are the biggest ways that FPSAC has changed since those earliest conferences you attended in the early 90s?}

\quot{F.B.}{What impresses me is the larger scope of the meeting, both mathematically and in the number of participants. There is no doubt that the papers submitted represent what is most exciting that is going on in the field on the international front. Young forefront researchers are clearly being selected for talks and posters. There is also a shift from more simple enumerative problems to deeper algebraic combinatorics questions and their natural ties to many subjects, such as representation theory, algebraic geometry, algebraic topology, theoretical physics, and more.}

\bio{H\'{e}l\`{e}ne Barcelo}{H\'{e}l\`{e}ne Barcelo is a professor emerita in the School of Mathematical and Statistical Sciences at Arizona State University, and from 2008 to 2026 was the Deputy Director of the Simons Laufer Mathematical Sciences Institute (SLMath, formerly MSRI). She received her PhD in math from the University of California, San Diego in 1988 under the supervision of Adriano Garsia. Her research is in algebraic and topological combinatorics, including representation theory, homotopy theory, and hyperplane arrangements.  She organized the 13th FPSAC in Tempe, Arizona in 2001, and was a member of the FPSAC executive committee in the 1990s and 2000s.}

\quot{S.H.}{How did you first get involved with FPSAC?}

\quot{H.B.}{I attended the conference on enumerative combinatorics in Montr\'{e}al in 1985 which was in many ways the precursor to FPSAC. I was a graduate student at the time, and when I heard that there was going to be a combinatorics conference in Montr\'{e}al, I forced my advisor, Adriano Garsia, to come to it. That conference was clearly the beginning of something; people clearly thought, ``we should repeat this.'' At the time, we did not think necessarily every year, but we thought we should have this combinatorics conference at a certain frequency. I was later definitely at the conference in Bordeaux in 1991 as well.}

\quot{S.H.}{And how many FPSAC conferences would you estimate you have attended?}

\quot{H.B.}{I really don't know. In the 1990s and 2000s, I attended many FPSACs in a row, at least ten years in a row. At some point in the 90s, I had attended so many FPSACs, someone told me they were trying to ``beat'' my attendance. After 2008, my attendance decreased significantly because I became deputy director at MSRI and my summers were very busy.}

\quot{S.H.}{You organized the 13th FPSAC in Tempe, Arizona in 2001. How did you decide to do that, and what was organizing like? Any memorable moments?}

\quot{H.B.}{In those days, na\"{i}vet\'{e} was much more present. Somebody, perhaps Maylis Delest, said, ``maybe having FPSAC in Arizona would be a good idea,'' and I said, ``oh yes!'' Then I started to inquire about what it would take at my institution and what kind of support I would have. At the time, the chair of the math department at ASU was Tom Trotter, a graph theorist, and he was supportive. But it was so much work! The other person who really helped me, not with the local organization but with the program and the proceedings, was Volkmar Welker. A lot of people came to that FPSAC, and I think the thing they will remember most about it was that it was bloody hot. We held the conference in May to try to avoid the worst of the summer, but there was still a heat wave that week.}

\quot{S.H.}{During the period when you were on the FPSAC executive committee, there was a requirement that all FPSAC submissions have their abstracts written in both English and French. (This requirement was removed in 2016.) Can you talk a little about that language requirement?}

\quot{H.B.}{We had to fight so much for that language requirement in the beginning, because of course people did not want to write their abstracts in two languages. You have to remember, when FPSAC started in the 80s, English was not quite the dominant language in science that it is today. The intent of that requirement, and allowing submissions in French, was to make sure that people from the French-speaking community were not at a disadvantage if they could not write well in English.}

\quot{S.H.}{How do you feel that FPSAC has evolved over the years?}

\quot{H.B.}{There are many young people now attending the conference, and not as many of the older people. It is good that so many young people are attending, but I also wish that the older people would attend more, like they did in the past, to mingle and talk with the younger people. In the past, pretty much every combinatorialist would attend FPSAC -- the France gang, the MIT gang, and so on.}

\quot{S.H.}{For many years, you served as deputy director at SLMath (formerly MSRI) in Berkeley. There you witnessed countless conferences, workshops, et cetera, be organized, in many different areas of math. How do you think FPSAC compares to these SLMath programs? Are there things that FPSAC could learn from them?}

\quot{H.B.}{It is very different to organize a conference at MSRI, which provides already all the organizational infrastructure to the mathematicians, so that they just have to concentrate on the scientific content. When you organize a FPSAC conference, you have to do everything yourself. It would be quite difficult and expensive, but if FPSAC could have some kind of permanent infrastructure, that would make things much easier for the organizers. At least, people applying to organize a FPSAC should make sure to get the support of their institutions even before they apply.}

\quot{S.H.}{Which FPSACs do you remember the most?}

\quot{H.B.}{I remember well the Bordeaux conference -- it was just a very warm event. I also strongly remember the conference in Taormina, Sicily. I helped Luisa Carini with the local organization of that conference, and it was quite an experience!}

\bio{Nantel Bergeron}{Nantel Bergeron is the York Research Chair in Applied Algebra and Distinguished Research Professor in the Department of Mathematics and Statistics of York University.  He received his PhD in math from the University of California, San Diego in 1990 under the supervision of Adriano Garsia. He has made important contributions to many topics in algebraic combinatorics, including Coxeter groups, Schubert calculus, Hopf algebras, quasi-symmetric functions, and Macdonald polynomials. He organized the 10th FPSAC in Toronto in 1998, and was a member of the FPSAC executive committee from 2002 to 2018.}

\quot{S.H.}{Your first FPSAC was the one in Bordeaux in 1991, which you attended together with your brother Fran\c{c}ois. Can you talk a little about that first FPSAC?}

\quot{N.B.}{Fran\c{c}ois and I had a joint contribution with Robert Howlett and Donald Taylor, entitled ``A decomposition of the descent algebra of finite Coxeter groups,'' and I also presented ``Une approche combinatoire des polyn\^{o}mes de Schubert.'' Looking back, those two titles already contain themes that would remain central to my work: Coxeter groups, descent algebras and Schubert calculus. \\ \indent Although Bordeaux was my first conference carrying the FPSAC name, it was not my first encounter with this community. As an undergraduate in the 1980s, I attended an earlier conference in Montr\'{e}al where I first met Adriano Garsia, who would become my advisor. Many of the same people were involved in organizing it, so for me Bordeaux felt like a continuation of that earlier Montr\'{e}al meeting. The community had found a name and begun building the conference series that became FPSAC. \\ \indent It was a relatively young and intimate conference series then, and many interesting things were happening. Mireille Bousquet-M\'{e}lou defended her thesis during the conference, almost as a parallel session. Her advisor, G\'{e}rard Viennot -- who later preferred to use the name Xavier Viennot -- had a great influence on me. He had a remarkable way of giving talks: with a simple picture, he could make a result seem both completely clear and almost obvious, without diminishing its depth.}

\quot{S.H.}{And how many FPSACs would you estimate you've attended since?}

\quot{N.B.}{From my first FPSAC in 1991 through the 2026 conference,
there were 36 editions. I was involved in one way or another in 33 of them. I actually attended 26, and in 7 additional years I contributed without being able to attend.}

\quot{S.H.}{How did you come to organize the 10th FPSAC in Toronto in 1998? What was the experience of organizing the conference like?}

\quot{N.B.}{At the time, there were only a small number of strong, visible groups in algebraic combinatorics: Bordeaux, Montr\'{e}al and Paris, together with groups in a few American cities such as San Diego, Boston and Minneapolis. FPSAC was then alternating between Europe and North America. In fact, at the very beginning there had even been talk of alternating simply between Bordeaux and Montr\'{e}al, but interest in the meeting grew much faster than those original plans. \\ \indent My experiences in Bordeaux and Montr\'{e}al made me realize how important a conference like FPSAC could be in establishing a research group in algebraic combinatorics. I began running a weekly seminar at York University and inviting renowned speakers. By then many cities wanted to organize FPSAC, so one had to think three or four years ahead. At the Paris meeting in 1995, I proposed Toronto for 1998, respecting the Europe--North America alternation. \\ \indent I had only recently obtained a permanent position, but I had all the pieces needed to make a strong proposal. The Fields Institute was the ideal institutional partner. Ian Goulden was nearby in Waterloo, much more established than I was, and was willing to help and to chair the program committee. Frank Sottile, with whom I had begun a wonderful collaboration, was extremely helpful in every aspect of the project. Between the three of us, we handled almost everything, from the final program selection -- where Ian's experience and guidance were invaluable -- to the banquet. At that time, a successful proposal had to demonstrate an ability to obtain funding, name a strong and established program committee chair, and present a convincing scientific plan, including whom one hoped to invite. \\ \indent We were very deliberate about the invited speakers. FPSAC was already an open-minded and diverse conference for its time, and we wanted diversity of research areas, gender and geography, particularly a balance between Europe and North America. We also wanted someone from outside algebraic combinatorics who used combinatorial ideas in another discipline. Bernard Derrida, coming from physics, was my ``outsider'' choice. The complete list was Georgia Benkart, Peter Cameron, Patrick Dehornoy, Bernard Derrida, Persi Diaconis, Chris Godsil, Ken Ono, Jean-Yves Thibon and Bernd Sturmfels. I remain very proud of that selection; several were already, or would become, international stars. \\ \indent I did choose to organize the banquet myself. Previous FPSAC banquets had always been special, and I did not want Toronto to miss that tradition. I rented Casa Loma for the evening, so participants could wander freely through the castle and gardens. I selected the food and wine and arranged a small chamber-music performance. It was tremendous fun, and I learned a surprising amount about organizing an event: venues, permits, legal requirements, catering and personnel.}

\quot{S.H.}{You were a member of the FPSAC executive committee (which had various previous names, like the ``executive office'' of the ``permanent program committee'') from 2002 to 2018. What caused you to resign in 2018?}

\quot{N.B.}{I had become concerned that FPSAC was not doing enough to encourage diversity or to make room for new people. Throughout my committee service, I had pushed the original philosophy that established and young mathematicians should meet as one community, and that the program should represent different people, institutions, countries and areas. By 2018, I felt that too many established mathematicians were giving contributed talks and that too much of the program came repeatedly from the same institutions. The resulting lack of diversity was contrary to what I believed FPSAC should be, and it eventually led me to resign. Those founding principles mattered to me, and, strangely, some of them need to be reaffirmed today.}

\quot{S.H.}{Speaking of promoting young researchers, starting with the conference in Vancouver in 2004 (where you were on both program and organizing committees), a ``best student paper'' award is given out each FPSAC. How did this award come to be, and has it been successful in its goals?}

\quot{N.B.}{As far as I remember, the award was proposed during the preparation of the Vancouver conference, although I cannot say with certainty who suggested it first. Jonathan Borwein and Stephanie van Willigenburg co-chaired the program committee, and Stephanie also chaired the graduate-prize committee. Lauren Williams was the first recipient. Each conference was expected to include a member of the permanent committee as a liaison; I was the liaison in Vancouver, and I immediately embraced the idea. The permanent committee adopted it and encouraged every subsequent conference to continue it. It was a wonderful idea. Looking at the list of winners now, one sees how many went on to become outstanding mathematicians. The prize did exactly what it was supposed to do: identify and encourage remarkable work at an early stage.}

\quot{S.H.}{Do you have a favorite FPSAC memory?}

\quot{N.B.}{It is difficult to choose just one. Bordeaux in 1991 has a special place because it was my first FPSAC: the warmth of the Bordeaux community and our visit to Viennot's vineyard remain vivid. The following year in Montr\'{e}al, Viennot's lecture accompanied by G\'{e}rard Duchamp on the violin showed me how clear, deep and beautiful mathematical communication could be. Toronto in 1998 is also unforgettable, especially seeing the community explore Casa Loma during the banquet we had worked so hard to create. Chile was influential in a different way: strengthening my connection
with the Colombian group had a lasting effect on both my research and my cultural perspective. \\ \indent But my favorite ``memory'' may really be something repeated: arriving in a new city each year, seeing old friends, meeting students who are attending for the first time, and watching some of them return later as established researchers and organizers. That continuity across generations is what FPSAC means to me.}

\bio{Christian Krattenthaler}{Christian Krattenthaler is a retired professor in the Faculty of Mathematics at the University of Vienna (Universit\"{a}t Wien). He received his PhD in math from the University of Vienna in 1983 under the supervision of Johann Cigler. He has made many fundamental contributions to combinatorial enumeration, with a focus on generating function techniques, hypergeometric series, and evaluations of determinants. He organized the 9th FPSAC in Vienna, Austria in 1997, and was a member of the FPSAC executive committee from 2003 to 2018.}

\quot{S.H.}{The first FPSAC you attended was Bordeaux 1991. How did you hear about this conference, and what motivated you to attended?}

\quot{C.K.}{A possibility, which is mostly speculation, is the following. At the time I was corresponding with Bruce Sagan concerning bijective problems that Richard Guy had raised. (Eventually, that became the paper with us three as co-authors.) During the correspondence we discovered that both of us are musicians, Bruce on the violin, myself on the piano. So Bruce told me that he will go to FPSAC in Bordeaux, and asked whether this would give the opportunity to play some pieces together. We did that during the conference week. But I can't tell whether it was indeed Bruce who made me aware of the conference, or whether I knew it through some other source.}

\quot{S.H.}{How many FPSAC conferences have you attended?}

\quot{C.K.}{I attended all of them except: the first two (as essentially everybody else); Paris 1995; Moscow 2000; Vancouver 2004; Ramat-Gan 2021 (I followed online); and Seattle 2026.}

\quot{S.H.}{Do you think that you are the person who has attended the most FPSACs?}

\quot{C.K.}{I am sure. In the 1990s, Tom Roby and Masao Ishikawa were ``in the lead.'' They later ``fell behind'': Tom because he went into education for some time, and Ishikawa's attendance became less regular in the past years.}

\quot{S.H.}{At the request of Daniel Krob (then the leading member of the FPSAC executive committee), you organized the 9th FPSAC in Vienna in 1997. Were there any unexpected challenges to organizing that conference?}

\quot{C.K.}{As far as I remember, there were no unexpected challenges. What was unexpected was the fact that the number of participants dropped from 200 two years prior to 130--150 that year. (Fortunately the hotels where we had blocked rooms did not complain.) It remained like that for several years, and it was only starting with the FPSAC in Linz in 2009 that numbers  jumped again to over 200.  As for an explanation for the decrease in attendance, I can again only speculate: in the beginning, FPSAC was the only big conference in the field; starting with the mid-1990s, many more conferences and meetings were organized. The subsequent rise of participant numbers may have its explanation in two things: first, the NSF/NSA funding for young participants from the US, which did not exist in 1997, for example; and second, I guess, the growth of our community.}

\quot{S.H.}{What are the biggest changes you've observed with FPSAC over the many years that you've been attending?}

\quot{C.K.}{There have been several big changes. Participant numbers. (I don't know why the numbers of participants are not on the official FPSAC website\footnote{\url{https://fpsac.org/}} anymore. The older website\footnote{\url{http://www-igm.univ-mlv.fr/~fpsac/}} shows the numbers up to the FPSAC in Reykjavik in 2011.) Number of submissions. (In Vienna there were 85 submissions.) Formalized procedures. (Overall it is good to have rules. Sometimes I get the impression that one can also overdo things.) The average quality of presentations has increased enormously when I compare, say, the talks in Vienna with the ones in recent FPSACs. (I mean both scientific quality and style of presentation.) \\ \indent Thematically, there has been certainly a shift towards more algebraic topics, at the loss of the ``formal power series'' part in the title of the conference. In my opinion, this development needs to be monitored. For me, FPSAC is so fascinating and fruitful since it brings people from very diverse research areas together who still have a common ground to talk to each other. If the mix of topics becomes too lop-sided, then FPSAC risks losing communities, with the result of narrowing down. This would really be a pity. In that regard, the program committee has a big responsibility, in particular the program chairs. This does not always work equally well every year. (A drastic example is ``Counting factorizations of Coxeter elements into products of reflections'' by Guillaume Chapuy and Christian Stump, which was submitted to and rejected from whatever FPSAC. I don't know an article that better fits the scope of FPSAC -- an enumerative problem motivated by algebra, solved by the skillful use of generating functions!) \\ \indent I find it also too bad that we have essentially lost the community interested in asymptotics. These asymptotics people are also very much interested in generating functions -- and, thus, power series, even though not formal in this context -- and therefore we have something to talk to each other about. \\ \indent I know that there are some who think that ``FPS'' should be removed from the title (and the conference). I can testify that already in 2003 at the FPSAC in Link\"{o}ping, Anders Bj\"{o}rner suggested to get rid of ``FPS.'' With Pierre Leroux being the leading figure at that time, this of course had no chance to get through. As I said, I disagree sharply with Anders on that point. It would be a pity to narrow the scope of the conference down.}

\quot{S.H.}{Do you have any favorite FPSAC memories?}

\quot{C.K.}{There are too many. Let me restrict myself to my first FPSAC, 1991 in Bordeaux. That was the second international conference that I was attending. I saw many famous people for the first time: G\'{e}rard Viennot (that was before he switched to his second first name), who delivered a performance as magician of transparencies; Doron Zeilberger, who delivered a show entitled ``Identities in search of identity''; Philippe Flajolet; Pierre Cartier; Gilbert Labelle as invited speakers. Julian West presented an intriguing conjecture in his contributed talk that triggered the development of equations with catalytic variables (``invented'' by Doron on the flight back). Very influential for me was the contributed talk by Dominique Gouyou-Beauchamps which had also an intriguing conjecture. I managed to prove it and this became my first \emph{Memoirs of the  AMS} publication, and many subsequent papers of mine are also heavily influenced by that work. \\ \indent Mireille Bousquet-Melou's thesis defense took place during the week (so-to-speak as an unofficial part of the conference). I attended it. It was about polyominoes, and -- obviously -- in French. Since I had no idea of French then, I did not understand much. It was nevertheless an interesting experience. \\ \indent The conference dinner took place in some castle, and Viennot invited several people to his house in some village outside of Bordeaux. At this occasion, Doron received a prize in the form of 10 bottles of Viennot wine for his ``Strahler number bijection.'' And, yes, Bruce Sagan and myself did some violin/piano music in a nearby building. At the FPSAC in Vienna 1997, we then did a joint recital.}

\bio{Mike Zabrocki}{Mike Zabrocki is a professor in the Department of Mathematics and Statistics at York University. He received his PhD in math from the University of California, San Diego in 1998 under the supervision of Adriano Garsia. His research is in algebraic combinatorics and representation theory, especially symmetric functions and their many variants. He was a co-organizer of the 18th FPSAC in San Diego in 2006, and since c.~2010 has been a member of the executive committee.}

\quot{S.H.}{What was the first FPSAC you attended? How did you become aware of it?}

\quot{M.Z.}{I first attended the FPSAC in 1999 in Barcelona. I suspect that I learned of FPSAC because other grad students were attending this conference.}

\quot{S.H.}{And how many FPSACs would you estimate that you've attended?}

\quot{M.Z.}{I went through the list of FPSACs, and it was easiest to remember the locations that I've attended. The answer is 21.}

\quot{S.H.}{Together with Jeff Remmel, you were one of the lead organizers of the FPSAC in San Diego in 2006. How did you decide to organize this FPSAC? Were there any major challenges in organization? Any memorable moments?}

\quot{M.Z.}{I had just gotten tenure at York University in Toronto and by the time I offered to organize I had attended 5 or so FPSACs. I saw it was an opportunity to invite people to come and experience San Diego.  Nantel Bergeron had organized a FPSAC in 1998 in Toronto, and he also went to grad school in San Diego and he was very supportive of making it possible for it to happen. I at least partially grew up in San Diego and I lived there 1988--1998 before moving to Canada. San Diego is a place that I really love to be, so I wanted to show that to the mathematics community. I had my first sabbatical in 2005--2006 so the timing was right for me to spend part of my sabbatical in San Diego and to raise my hand and see if I could make it happen. \\ \indent I approached Jeff Remmel because we needed the help of the university for logistics and administrative support and he was very good at pulling strings. The first challenge was me being a ``local organizer'' and only being a visitor in the department of mathematics. I had to ask for administrative help from the department, but not being someone who had a say in the department made it awkward for me to ask them to contribute. Most of our organizing committee were also ``not local.'' \\ \indent The second challenge we had was the infrastructure at UCSD.  There were no hotels close to the university at the time and there were restrictions on the dorms on campus that made it difficult to hold the event on campus. The university is far from the surrounding neighborhoods.  We decided to hold it in a hotel in Pacific Beach. The conference was much smaller back then, so I wouldn't even attempt to organize it that way now. Even then it was not clear how to manage to hold a conference of that size in San Diego. \\ \indent Most memorable moment for me: the opening reception was on the beach and my parents helped put out a huge spread with tables and decorations.}

\quot{S.H.}{How has FPSAC changed over the years that you've been attending?}

\quot{M.Z.}{The local organization has really upped the game every single year. The extra flourishes that the organizers add to the conference make the event so nice to attend and experience the location every year. But the growth of the conference has made the conference challenging to put on both logistically and because we have a low acceptance rate.}

\quot{S.H.}{Do you have any favorite FPSAC memories?}

\quot{M.Z.}{I look back at the list of the conferences that I have attended, and the answer if I have favorite FPSAC memories is ``yes.'' Usually it's one memorable feature from each one that I have attended: the banquet at Barcelona, the swimming pool in Tempe, the coffee and wine in Melbourne, et cetera. I think that we are lucky to be able to go to some really spectacular locations with this conference and I appreciate how the local organizers like to show attendees their ``home town.''}

\quot{S.H.}{What are the biggest challenges that FPSAC faces, and how is the executive committee trying to address these challenges?}

\quot{M.Z.}{The biggest problem FPSAC has now is that the conference is outgrowing its own structure. The number of abstracts submitted each year poses two challenges: reviewing them properly and making room in the conference for the great work being submitted. The fact that FPSAC doesn't have parallel sessions means we roughly have 27 talks and (depending on infrastructure) 80--100 posters. With current number submissions, there has been less than fifty percent acceptance rate for the last many years. The program committee is forced to cut great material. \\ \indent There needs to be a balance between showcasing younger mathematicians (who are there to network and interact with the community) and more established mathematicians (who will meet and interact with the younger generation). The executive committee voted a few years ago to decrease the number of invited speakers in hopes that we can insert more participation into the schedule. The change will trickle in this year and we will have to see if it helps. \\ \indent When the number of attendees is larger than 300, not all universities have appropriate rooms. We have to work hard to find local organizers who will make the conference happen. It always surprises me that every year we have groups willing to step up and bring the community to yet another new and surprising place.}

\bio{Bridget Tenner}{Bridget Tenner is a professor in the Department of Mathematics at DePaul University. She received her PhD in math from the Massachusetts Institute of Technology (MIT) in 2006 under the supervision of Richard Stanley. Her research is in enumerative and algebraic combinatorics, with a focus on permutations, especially permutation patterns. She organized the 26th FPSAC in Chicago in 2014, and since 2015 she has been a member of the FPSAC executive committee.}

\quot{S.H.}{What was the first FPSAC conference you attended? How did you become aware of it, and what motivated you to attend?}

\quot{B.T.}{2003 Link\"{o}ping University, Sweden. I was in grad school at the time, and my advisor Richard Stanley told me about it. I was shocked -- I might get to go to Sweden? and learn math? and have funding support to do so? My stepfather was an anthropologist and archaeologist, and I grew up familiar with the lifestyle of academia, but I hadn't yet put together that I might get to live that lifestyle myself. It was very exciting!}

\quot{S.H.}{You organized the 26th FPSAC in Chicago in 2014. (Incidentally, this was the first FPSAC I attended, and it certainly set the bar -- I will always remember the Goose Island conference dinner!) How did you decide to organize that FPSAC?}

\quot{B.T.}{Mike Zabrocki had approached me about organizing a Chicago FPSAC during FPSAC in Iceland a few years before. It was flattering to be asked and trusted -- both with pulling it off, and with showcasing my spin on conferences and combinatorics. I was also keen to further establish Chicago as a place for algebraic (and enumerative) combinatorics. There are a number of good universities in town, but the community had only recently been building up steam. I could see that we had some good things going here, and I was eager for other people to see them too. Plus, Chicago is great!}

\quot{S.H.}{How were the plenary speakers selected for the Chicago FPSAC? I'm especially interested in the choice of June Huh, who had then just finished his PhD.}

\quot{B.T.}{If I remember correctly, the plenary speakers were mostly worked out through the program committee chairs (Lou Billera and Bella Novik), with suggestions here and there from me and others. We worked to have a good balance along all kinds of dimensions. As we collected a list of names, they were sometimes partitioned by, for example, research area, so that we would be trying to find a person in probabilistic combinatorics, and we'd then narrow that list down to three people. \\ \indent The discussion where June Huh's name came up went, approximately, as follows. Lou: To the list of ``particularly interesting'' recent work, let me add June Huh, who just proved that the first kind Whitney numbers of realizable matroids are log concave, settling a very old conjecture of Rota, Welsh et al.~in the realizable case. Bella: I agree that June's work is a real breakthrough -- we should definitely have him on our list!}

\quot{S.H.}{At that FPSAC in Chicago, you recorded the talks and subsequently made the recordings available online. More recently, there has also been livestreaming of FPSAC talks. It can be tricky to balance the obvious benefits of allowing online participation versus the unique experience of an in person conference. What are your thoughts about that?}

\quot{B.T.}{This is a challenging thing. I prefer in-person events and participation whenever possible. They allow for unstructured interactions, which can sometimes be the most fruitful. I think there's also something meaningful about the energy in a room during a talk. And as a speaker, I can give a much more interactive -- and, I think, a stronger -- talk if there are people in the room with me. Plus I'm a big one for moving around and gesturing, and this is hard to do remotely! That being said, there are a lot of reasons why in-person participation might not be possible, and so a streaming option is important. I have yet to see hybrid events (talks, conferences, classes) go smoothly, and so my own preference is that when in-person and streaming are happening simultaneously, there should probably be a designated person in the lecture hall monitoring the stream, acting as the voice of any relevant questions or comments coming through the chat.}

\quot{S.H.}{What are the biggest changes to the conference you've observed in the years you've been attending FPSAC?}

\quot{B.T.}{It is growing a lot. The place I think about this the most is in terms of the submissions. A lot of great work gets rejected, simply because of time constraints. It is very unfortunate and can be hard on the submitters. Also on program committee members who have a huge load to lift.}

\quot{S.H.}{Do you have any favorite FPSAC memories?}

\quot{B.T.}{Long post-conference dinners with friends that I hadn't seen since the last conference! When it was in Chicago, some folks came over to my house after the last day and cooked me dinner. It was a welcome treat, and a lovely celebration together of hard work being done.}

\quot{S.H.}{How do you see FPSAC developing in the future? What challenges does the conference face, and how is the executive committee attempting to address them?}

\quot{B.T.}{I think the size and interest is going to be a big challenge. The no-parallel-talks nature of FPSAC is at the core of the conference's design, and yet it is a tough knock to get rejected year after year. The executive committee has made some attempts to brainstorm how to increase the number of presenters, but it can be very site-specific. Also, navigating in-person versus online versus hybrid participation is a big question mark, too. It'll be interesting to see how the role of AI impacts the conference, just as this is evolving in the broader mathematical community. Of course, a lot of these issues are in the hands of the local organizing and program committees, but they're important for all of us to think about.}

\bio{Jang Soo Kim}{Jang Soo Kim is a professor in the Department of Mathematics at Sungkyunkwan University. He received his PhD in math from Korea Advanced Institute of Science and Technology (KAIST) in 2009 under the supervision of Dongsu Kim. His research is in enumerative and algebraic combinatorics, studying things like tableaux, posets, polytopes, and orthogonal polynomials. He helped organize the 27th FPSAC in Daejeon, South Korea in 2015, and since 2019 he has been a member of the FPSAC executive committee.}

\quot{S.H}{What was the first FPSAC conference you attended? How did you become aware of it, and what motivated you to attend?}

\quot{J.S.K.}{My first FPSAC conference was FPSAC 2007 in Tianjin, China. I heard about FPSAC from Heesung Shin and Seunghyun Seo, who had attended an earlier FPSAC. My advisor, Dongsu Kim, recommended that I attend FPSAC 2007.}

\quot{S.H.}{And how many FPSACs have you attended since then?}

\quot{J.S.K.}{I've attended every FPSAC since then, although I participated online in FPSAC 2021, which was held in Israel. If online attendance counts, then I've never missed an FPSAC since 2007.}

\quot{S.H.}{You helped organize the 27th FPSAC in Daejeon, South Korea in 2015. How did you get involved with that? Were there unexpected challenges in organization?}

\quot{J.S.K.}{The main organizer of FPSAC 2015 was Dongsu Kim. Since he was my advisor and I had already attended FPSAC many times, it was quite natural for me to join the organizing committee that year. In the summer of 2015, there was an outbreak of MERS, so we were concerned about whether we could hold FPSAC as scheduled. Fortunately, the outbreak was not as serious as we had feared.}

\quot{S.H.}{FPSAC now regularly has many participants from North America, Europe, and Asia. Do you think it does a good job accommodating this wide geographic diversity? Is that something the executive committee is focused on?}

\quot{J.S.K.}{Over the last 15 years, FPSAC venues have been distributed evenly across North America, Europe, and Asia, mostly in different countries. Based on these statistics, I think the executive committee has done well in recent years with regard to geographic diversity. It is part of our mission to find places where FPSAC has not been held before. Doing so will make the community more diverse and attract people from countries where research in the areas covered by FPSAC is less developed.}

\quot{S.H.}{Do you have a favorite FPSAC memory?}

\quot{J.S.K.}{The banquet at FPSAC 2018 in Dartmouth was fantastic! They provided lobsters and an unlimited supply of alcohol. We all danced together to music played by a band. I think dancing at FPSAC banquets has since become a tradition.}

\quot{S.H.}{In 2024, you started the FPSAC archive,\footnote{\url{https://fpsac-archive.github.io/}} which collects detailed information about all the FPSAC conferences. What motivated you to do this? Have you learned anything interesting about the history of FPSAC from this archive project?}

\quot{J.S.K.}{In December 2023, someone asked the executive committee whether papers from past conferences were available. At the time, many papers were missing. I thought it would be useful to make all FPSAC papers available online, so I decided to collect the missing ones. I contacted former executive committee members, and Christophe Reutenauer, Volker Strehl, and Christian Krattenthaler sent me some of the missing documents. Christian Krattenthaler helped me greatly by scanning his hard copies of previous FPSAC proceedings. The proceedings are quite useful for seeing how FPSAC has evolved. We can see how the numbers of papers, and how they were published, have changed over time. From 1988 to 2002, each FPSAC had proceedings containing selected papers presented at the conference, sometimes accompanied by original papers. From 2008 to 2016, the FPSAC proceedings were published in special issues of DMTCS (\emph{Discrete Mathematics \& Theoretical Computer Science}). Since 2017, the proceedings have been published in SLC (\emph{S\'{e}minaire Lotharingien de Combinatoire}).}

\quot{S.H.}{How do you see FPSAC developing in the future? What challenges does the conference face, and how is the executive committee attempting to address them?}

\quot{J.S.K.}{The conference is getting bigger and bigger. For the 2026 conference, where I was a program committee co-chair, we received 252 submissions. I expect this number to grow even more rapidly because of the availability of AI tools. The executive committee is exploring ways to address this growth.}

\bibliographystyle{abbrv}
\bibliography{main}

\end{document}